\documentstyle[11pt,bezier]{article} 
\newcommand{\qbezier}{\bezier{120}}
\def\DottedCircle{
\bezier{4}(0.966,-0.259)(1.04,0)(0.966,0.259)
\bezier{4}(0.966,0.259)(0.897,0.518)(0.707,0.707)
\bezier{4}(0.707,0.707)(0.518,0.897)(0.259,0.966)
\bezier{4}(0.259,0.966)(0,1.04)(-0.259,0.966)
\bezier{4}(-0.259,0.966)(-0.518,0.897)(-0.707,0.707)
\bezier{4}(-0.707,0.707)(-0.897,0.518)(-0.966,0.259)
\bezier{4}(-0.966,0.259)(-1.04,0)(-0.966,-0.259)
\bezier{4}(-0.966,-0.259)(-0.897,-0.518)(-0.707,-0.707)
\bezier{4}(-0.707,-0.707)(-0.518,-0.897)(-0.259,-0.966)
\bezier{4}(-0.259,-0.966)(0,-1.04)(0.259,-0.966)
\bezier{4}(0.259,-0.966)(0.518,-0.897)(0.707,-0.707)
\bezier{4}(0.707,-0.707)(0.897,-0.518)(0.966,-0.259)
}
\def\Endpoint[#1]{
\ifcase#1
\put(1,0){\circle*{0.15}}
\or\put(0.866,0.5){\circle*{0.15}}
\or\put(0.5,0.866){\circle*{0.15}}
\or\put(0,1){\circle*{0.15}}
\or\put(-0.5,0.866){\circle*{0.15}}
\or\put(-0.866,0.5){\circle*{0.15}}
\or\put(-1,0){\circle*{0.15}}
\or\put(-0.866,-0.5){\circle*{0.15}}
\or\put(-0.5,-0.866){\circle*{0.15}}
\or\put(0,-1){\circle*{0.15}}
\or\put(0.5,-0.866){\circle*{0.15}}
\or\put(0.866,-0.5){\circle*{0.15}}
\fi}
\def\Arc[#1]{
\thicklines                     % this can be changed!
\ifcase#1
\bezier{25}(0.966,-0.259)(1.04,0)(0.966,0.259)
\or
\bezier{25}(0.966,0.259)(0.897,0.518)(0.707,0.707)
\or
\bezier{25}(0.707,0.707)(0.518,0.897)(0.259,0.966)
\or
\bezier{25}(0.259,0.966)(0,1.04)(-0.259,0.966)
\or
\bezier{25}(-0.259,0.966)(-0.518,0.897)(-0.707,0.707)
\or
\bezier{25}(-0.707,0.707)(-0.897,0.518)(-0.966,0.259)
\or
\bezier{25}(-0.966,0.259)(-1.04,0)(-0.966,-0.259)
\or
\bezier{25}(-0.966,-0.259)(-0.897,-0.518)(-0.707,-0.707)
\or
\bezier{25}(-0.707,-0.707)(-0.518,-0.897)(-0.259,-0.966)
\or
\bezier{25}(-0.259,-0.966)(0,-1.04)(0.259,-0.966)
\or
\bezier{25}(0.259,-0.966)(0.518,-0.897)(0.707,-0.707)
\or
\bezier{25}(0.707,-0.707)(0.897,-0.518)(0.966,-0.259)
\fi}
\def\DottedArc[#1]{
\ifcase#1
\bezier{4}(0.966,-0.259)(1.04,0)(0.966,0.259)
\or
\bezier{4}(0.966,0.259)(0.897,0.518)(0.707,0.707)
\or
\bezier{4}(0.707,0.707)(0.518,0.897)(0.259,0.966)
\or
\bezier{4}(0.259,0.966)(0,1.04)(-0.259,0.966)
\or
\bezier{4}(-0.259,0.966)(-0.518,0.897)(-0.707,0.707)
\or
\bezier{4}(-0.707,0.707)(-0.897,0.518)(-0.966,0.259)
\or
\bezier{4}(-0.966,0.259)(-1.04,0)(-0.966,-0.259)
\or
\bezier{4}(-0.966,-0.259)(-0.897,-0.518)(-0.707,-0.707)
\or
\bezier{4}(-0.707,-0.707)(-0.518,-0.897)(-0.259,-0.966)
\or
\bezier{4}(-0.259,-0.966)(0,-1.04)(0.259,-0.966)
\or
\bezier{4}(0.259,-0.966)(0.518,-0.897)(0.707,-0.707)
\or
\bezier{4}(0.707,-0.707)(0.897,-0.518)(0.966,-0.259)
\fi}
\def\Chord[#1,#2]{
\thinlines
\ifnum#1>#2\Chord[#2,#1]
\else\ifnum#1<#2
\ifcase#1
\ifcase#2
\or\qbezier(1,0)(0.516,0.138)(0.866,0.5)
\or\qbezier(1,0)(0.45,0.26)(0.5,0.866)
\or\qbezier(1,0)(0.327,0.327)(0,1)
\or\qbezier(1,0)(0.179,0.311)(-0.5,0.866)
\or\qbezier(1,0)(0.0536,0.2)(-0.866,0.5)
\or\put(1, 0){\line(-2, 0){2}}
\or\qbezier(1,0)(0.0536,-0.2)(-0.866,-0.5)
\or\qbezier(1,0)(0.179,-0.311)(-0.5,-0.866)
\or\qbezier(1,0)(0.327,-0.327)(0,-1)
\or\qbezier(1,0)(0.45,-0.26)(0.5,-0.866)
\or\qbezier(1,0)(0.516,-0.138)(0.866,-0.5)
\fi
\or\ifcase#2\or
\or\qbezier(0.866,0.5)(0.378,0.378)(0.5,0.866)
\or\qbezier(0.866,0.5)(0.26,0.45)(0,1)
\or\qbezier(0.866,0.5)(0.12,0.446)(-0.5,0.866)
\or\qbezier(0.866,0.5)(0,0.359)(-0.866,0.5)
\or\qbezier(0.866,0.5)(-0.0536,0.2)(-1,0)
\or\put(0.866, 0.5){\line(-5, -3){1.73}}
\or\qbezier(0.866,0.5)(0.146,-0.146)(-0.5,-0.866)
\or\qbezier(0.866,0.5)(0.311,-0.179)(0,-1)
\or\qbezier(0.866,0.5)(0.446,-0.12)(0.5,-0.866)
\or\qbezier(0.866,0.5)(0.52,0)(0.866,-0.5)
\fi
\or\ifcase#2\or\or
\or\qbezier(0.5,0.866)(0.138,0.516)(0,1)
\or\qbezier(0.5,0.866)(0,0.52)(-0.5,0.866)
\or\qbezier(0.5,0.866)(-0.12,0.446)(-0.866,0.5)
\or\qbezier(0.5,0.866)(-0.179,0.311)(-1,0)
\or\qbezier(0.5,0.866)(-0.146,0.146)(-0.866,-0.5)
\or\put(0.5, 0.866){\line(-3, -5){1}}
\or\qbezier(0.5,0.866)(0.2,-0.0536)(0,-1)
\or\qbezier(0.5,0.866)(0.359,0)(0.5,-0.866)
\or\qbezier(0.5,0.866)(0.446,0.12)(0.866,-0.5)
\fi
\or\ifcase#2\or\or\or
\or\qbezier(0,1.)(-0.138,0.516)(-0.5,0.866)
\or\qbezier(0,1.)(-0.26,0.45)(-0.866,0.5)
\or\qbezier(0,1.)(-0.327,0.327)(-1,0)
\or\qbezier(0,1.)(-0.311,0.179)(-0.866,-0.5)
\or\qbezier(0,1.)(-0.2,0.0536)(-0.5,-0.866)
\or\put(0, 1){\line(0, -2){2}}
\or\qbezier(0,1.)(0.2,0.0536)(0.5,-0.866)
\or\qbezier(0,1.)(0.311,0.179)(0.866,-0.5)
\fi
\or\ifcase#2\or\or\or\or
\or\qbezier(-0.5,0.866)(-0.378,0.378)(-0.866,0.5)
\or\qbezier(-0.5,0.866)(-0.45,0.26)(-1,0)
\or\qbezier(-0.5,0.866)(-0.446,0.12)(-0.866,-0.5)
\or\qbezier(-0.5,0.866)(-0.359,0)(-0.5,-0.866)
\or\qbezier(-0.5,0.866)(-0.2,-0.0536)(0,-1)
\or\put(-0.5, 0.866){\line(3, -5){1}}
\or\qbezier(-0.5,0.866)(0.146,0.146)(0.866,-0.5)
\fi
\or\ifcase#2\or\or\or\or\or
\or\qbezier(-0.866,0.5)(-0.516,0.138)(-1,0)
\or\qbezier(-0.866,0.5)(-0.52,0)(-0.866,-0.5)
\or\qbezier(-0.866,0.5)(-0.446,-0.12)(-0.5,-0.866)
\or\qbezier(-0.866,0.5)(-0.311,-0.179)(0,-1)
\or\qbezier(-0.866,0.5)(-0.146,-0.146)(0.5,-0.866)
\or\put(-0.866, 0.5){\line(5, -3){1.73}}
\fi
\or\ifcase#2\or\or\or\or\or\or
\or\qbezier(-1,0)(-0.516,-0.138)(-0.866,-0.5)
\or\qbezier(-1,0)(-0.45,-0.26)(-0.5,-0.866)
\or\qbezier(-1,0)(-0.327,-0.327)(0,-1)
\or\qbezier(-1,0)(-0.179,-0.311)(0.5,-0.866)
\or\qbezier(-1,0)(-0.0536,-0.2)(0.866,-0.5)
\fi
\or\ifcase#2\or\or\or\or\or\or\or
\or\qbezier(-0.866,-0.5)(-0.378,-0.378)(-0.5,-0.866)
\or\qbezier(-0.866,-0.5)(-0.26,-0.45)(0,-1)
\or\qbezier(-0.866,-0.5)(-0.12,-0.446)(0.5,-0.866)
\or\qbezier(-0.866,-0.5)(0,-0.359)(0.866,-0.5)
\fi
\or\ifcase#2\or\or\or\or\or\or\or\or
\or\qbezier(-0.5,-0.866)(-0.138,-0.516)(0,-1)
\or\qbezier(-0.5,-0.866)(0,-0.52)(0.5,-0.866)
\or\qbezier(-0.5,-0.866)(0.12,-0.446)(0.866,-0.5)
\fi
\or\ifcase#2\or\or\or\or\or\or\or\or\or
\or\qbezier(0,-1.)(0.138,-0.516)(0.5,-0.866)
\or\qbezier(0,-1.)(0.26,-0.45)(0.866,-0.5)
\fi
\or\ifcase#2\or\or\or\or\or\or\or\or\or\or
\or\qbezier(0.5,-0.866)(0.378,-0.378)(0.866,-0.5)
\fi\fi\fi\fi}
\def\FullChord[#1,#2]{
\Endpoint[#1]
\Endpoint[#2]
\Arc[#1]
\Arc[#2]
\Chord[#1,#2]
}
\def\EndChord[#1,#2]{
\Endpoint[#1]
\Endpoint[#2]
\Chord[#1,#2]
}
\def\Picture#1{
\begin{picture}(2,1)(-1,-0.167)
#1
\end{picture}
}
\def\DottedChordDiagram[#1,#2]{
\Picture{\DottedCircle \FullChord[#1,#2]}
}
\newcommand{\e}{\varepsilon} 
\newtheorem{theor}{Theorem}[section] 
\newtheorem{thm}[theor]{Theorem} 
\newtheorem{lemma}[theor]{Lemma} 
\newtheorem{prop}[theor]{Proposition} 
\newtheorem{corol}[theor]{Corollary} 
       \newtheorem{rem}[theor]{Remark} 
\newcommand{\begproof}{\noindent \\ {\em Proof\/.} } 
\newcommand{\eproof}{ {\quad \ \hfill $\Box$ \\ \medskip}} 
\renewcommand{\a}{\alpha} 
\renewcommand{\b}{\beta} 
\newcommand{\g}{{\gamma}} 
\renewcommand{\d}{\delta}
\newcommand{\s}{\sigma} 
\newcommand{\be}{\begin{equation}} 
\newcommand{\ee}{\end{equation}} 
\newcommand{\bea}{\begin{eqnarray}} 
\newcommand{\eea}{\end{eqnarray}} 
\newcommand{\bean}{\begin{eqnarray*}} 
\newcommand{\eean}{\end{eqnarray*}} 
\newcommand{\text}[1]{ {\rm {#1} } }           % only for LaTeX 2.09 
\newcommand{\Bbb}{\bf}                       % only for LaTeX  2.09  
\newcommand{\mathbf}{\bf}                    % only for LaTeX  2.09  
\newcommand{\mathrm}{\rm}                    % only for LaTeX  2.09  
 
\newcommand{\C}{{\Bbb C}} 
\newcommand{\R}{{\Bbb R}}

\newcommand{\D}{\Delta}
 
\newcommand{\Th}{\Theta}
\newcommand{\bcA}{\bar{\cA}}
\newcommand{\cA}{{\cal A}} 
\newcommand{\cB}{{\cal B}} 
 
\newcommand{\cF}{{\cal F}} 
 
\newcommand{\cD}{{\cal D}}

\newcommand{\cP}{{\cal P}}

\newcommand{\cU}{{\cal U}} 
\newcommand{\cV}{{\cal V}} 
\newcommand{\cW}{{\cal W}} 
\newcommand{\Z}{{\Bbb Z}} 

\newcommand{\BN}{Bar-Natan} 
\newcommand{\Gar}{Garoufalidis}
\newcommand{\AC}{Alexander-Conway}
\newcommand{\Fn}{Feynman} 
\newcommand{\Fd}{\Fn\ diagram} 
\newcommand{\gl}{\mbox{$gl(1|1)$}} 
\newcommand{\hW}{\widehat{W}}
\newcommand{\ra}{\rightarrow} 
 
\newcommand{\Vas}{Vassiliev} 
\newcommand{\VI}{\Vas\ invariants} 
\newcommand{\ki}{knot invariant} 
\newcommand{\ws}{weight system} 

\newcommand{\half}{\frac{1}{2}}
\newcommand{\ov}{\overline}

\newcommand{\MM}{Melvin-Morton}

\begin{document} 

\title{
\bf \MM\ conjecture and primitive \Fd s}

\author{
Arkady Vaintrob
\\ {\em Department of Mathematics, University of Utah%
\thanks{on leave from New Mexico State University, Las
Cruces, NM 88003}%
,} 
\\ {\em Salt Lake City, UT 84112, USA \qquad e-mail: {\tt
vaintrob@math.utah.edu} }
}

%\date{\today}
\date{}
\maketitle

\begin{abstract}  
We give a very short proof of the \MM\ conjecture relating
the colored Jones polynomial and the Alexander polynomial of knots.
The proof is based on explicit evaluation of
the corresponding weight systems on primitive elements
of the Hopf algebra of chord diagrams which, in turn, follows
from simple identities between four-valent tensors  on
the Lie algebra $sl_2$ and the Lie superalgebra $\gl$.
This shows that the miraculous connection between the Jones and
Alexander invariants follows from the similarity
(supersymmetry) between $sl_2$ and $\gl$.
\end{abstract} 
 
%\maketitle                %  only for LaTeX2e !!!

\section{Introduction}
It is well known that the Jones polynomial $V$ of knots can be constructed
by working with the standard two-dimensional
representation of the quantum group $sl_2$.
An arbitrary, $d$-dimen\-is\-onal, representation of $sl_{2,q}$
gives rise to a generalization of $V$, 
the so-called {\em colored Jones invariant\/} $V^d$
(see \cite{KR,RT} for more details).
Morton and Strickland \cite{MS} proved
that $V^d$ is determined by all
the cablings of the Jones polynomial thus showing that
$V^d(K,t)$ can be computed combinatorially
by formulas generalizing the famous skein relation for $V(K,t)$.

The \AC\ polynomial
is another famous knot invariant that can be computed by using
a skein relation.
In \cite{MM} Melvin and Morton  studied pieces of $V^d(K,t)$ appearing
in a certain power series expansions and conjectured
that the Alexander polynomial of a knot is 
determined by its colored Jones invariant.
They proved that the coefficients $a_{in}$
of the power series expansion 
\be
V^d(K,e^z) = \sum_{i,n \ge 0} a_{in}(K)d^iz^n
\label{eq:expansion}
\ee
of $V^d$ in variables 
$d$ and $z=\log t$
are \Vas\ \ki s
of order $\le n$ and that $a_{in}=0$ if $i>2n$.
They formulated the following {\em \MM\ conjecture\/}
(later Morton~\cite{M} proved it for torus knots).  

\begin{thm} \label{thm:mm} \ \nopagebreak

(i) The coefficient ``matrix'' $(a_{in})$ of the expansion
(\ref{eq:expansion}) is upper triangular, i.e.
\be
 a_{in} = 0 \quad \text{for} \ i>n 
\label{eq:MM0}
\ee

(ii)  the leading (``diagonal'') term 
\be
V_0(K,z) = \sum_{n\ge 0} a_{nn}(K) z^n 
\label{eq:V0}
\ee
is the inverse of the renormalized Alexander polynomial of $K$, i.e.
\be
V_0(K,z) \cdot \frac{z}{e^{z/2} - e^{-z/2}} \Delta(K,e^z) = 1~.
\label{eq:MM}
\ee
\end{thm}

Rozansky \cite{Roz} derived relations~(\ref{eq:MM0}) and~(\ref{eq:MM})
from Witten's~\cite{W} path integral interpretation of
the Jones invariant  using quantum field theory tools.
\BN\ and  \Gar\ \cite{BNG} showed that it is enough to
 establish~(\ref{eq:MM0}) and~(\ref{eq:MM}) on the level of
 \Vas\ \ws  s and  proved the \MM\ conjecture by
finding combinatorial expressions for the \ws s corresponding to
the leading terms $a_{nn}$ of the expansion~(\ref{eq:expansion}) and the 
coefficients of the \AC\ polynomial.  
In~\cite{Krick} A.~Kricker
gave another proof based on the characterization~\cite{AK} of the 
subspace generated by the
\AC\ \ws s in terms of the cabling operations.
In~\cite{Vmm} we gave a different, simpler proof using
the recursion relations from~\cite{CV} and~\cite{FKV}
 for computing the Jones and the \AC\ \ws s. 

All these proofs are to various extent based on some miraculous
cancellations or coincidences occurring in computations of the Jones and
the \AC\ \ws s. This leaves one puzzled wondering whether
Theorem~\ref{thm:mm} could have been discovered without the
prior experimental work, wishing for a more conceptual 
explanation of the miracles. 

In this note we attempt to address this question. We give a very short
proof of the identities~(\ref{eq:MM0}) and~(\ref{eq:MM}) on the level
of \ws s that does not use any involved computations. The  idea is to
restrict everything  to primitive elements 
of the Hopf algebra of chord diagrams and to use the fact that the
colored Jones and the \AC\ invariants are the universal \VI\ coming
from the Lie algebra $sl_2$ and the Lie superalgebra $\gl$
(see~\cite{FKV}).  Among other things this allows us to ignore the
difference between the framed and unframed versions of the colored
Jones invariant. 
All that is needed to compute the  corresponding \ws s
on the primitive subspace are 
simple identities~(\ref{eq:lagrange}) and~(\ref{eq:ACiv})
 between invariant tensors of order four on
$sl_2$ and $\gl$.%
\footnote{In~\cite{CV}
and~\cite{FKV} these identities were used to derive recursive
formulas for the universal \ws s with values in $\cU(sl_2)$
and $\cU(\gl)$ respectively. }{}
 The superalgebra $\gl$ is the fermionic counterpart
of $sl_2$ (or rather $gl_2 = sl_2 \oplus \C$) which explains the
similarity between 
relations~(\ref{eq:lagrange}) and~(\ref{eq:ACiv}). However,
algebraically  $\gl$ is much simpler than $sl_2$, it is solvable and 
not {\em simple\/} as a Lie superalgebra. This is precisely the reason
why the Alexander invariant is weaker than the colored Jones invariant
being equivalent to its leading part~(\ref{eq:V0}).

The outline of the paper is as follows. 
In Section~2 we review the basic facts on
\VI\ and the Hopf algebra $\cA$
of chord diagrams. In Section~3 we establish some facts on
primitive elements of the algebra $\cA$. A construction assigning to
a Lie superalgebra with an invariant inner product a series of \Vas\
\ws s is recalled in Section~4.
In Section~5 the \MM\ conjecture is reduced to two statements
on the values of 
the Jones and the \AC\ \ws\ on the primitive elements of $\cA$
which are proved in Sections~6 and~7 respectively.

\section{\VI\ and \ws s}
All necessary information on \VI, chord diagrams  and
\ws s can be found, for example, in~\cite{BN,Konts,BL}.
The purpose of this section is to establish notation and
introduce some lesser standard terminology.

Let $\cV_n$ denote the space of \Vas\ \ki s of order $\le n$ with
values in a field $k$ (we usually assume that $k=\C$).
The space $\cV=\bigcup \cV_n$ of all
 \VI\/  is a filtered Hopf algebra with
the comultiplication induced by the operation of connected
sum of knots and the antipode given by the orientation
reversal. 

A {\em chord diagram\/} of order $n$ is an oriented circle 
(called the Wilson line) with $n$ disjoint pairs of points 
({\em chords\/}) up to an orientation
preserving diffeomorphism of the circle.  
Denote by $\cA_n$ the space generated by the  set  $\cD_n$ of
all chord diagrams with $n$ chords
modulo the following {\em four-term relations\/} 
\begin{equation}
% I\left(
\Picture{\DottedCircle\FullChord[1,8]\Arc[2]\FullChord[5,9]}
%\right) 
-
% I\left(
\Picture{\DottedCircle\FullChord[1,9]\Arc[2]\FullChord[5,8]}
% \right) 
+
%I\left(
\Picture{\DottedCircle\FullChord[2,5]\Arc[8]\FullChord[1,9]}
% \right) 
-
% I\left(
\Picture{\DottedCircle\FullChord[1,5]\Arc[8]\FullChord[2,9]}
% \right) 
= 0. \label{eq:4term}
\end{equation}\\[12pt]
The  space $\langle \cD_n \rangle$
modulo (\ref{eq:4term}) and the following {\em one-term relations\/}
\begin{equation}
\Picture{\DottedCircle\FullChord[2,10]\Arc[1]\Arc[0]\Arc[11]}
= 0\label{eq:1term}
\end{equation}\\[8pt]
will be denoted by $\bcA_n$.

The spaces 
$$\cA=\bigoplus_{n\ge 0}\cA_n \text{\ \ and\  \ }
\bcA=\bigoplus_{n\ge 0}\bcA_n
$$
have structures of  $\Z$-graded
commutative and co-commutative Hopf algebras
with the product given by the operation of the
connected sum of diagrams
(which is well-defined modulo four-term relations)
and the coproduct $\D: \cA \to \cA\otimes \cA$ by the formula
$$
\D(D) = \sum_{E\subset D} E \otimes (D\setminus E).
$$

The  projection $\cA \to \bcA$  splits and $\bcA$ 
can be embedded into $\cA$ as a Hopf subalgebra so that
$$\cA = \bcA[\Th],
$$ 
where $\Th$ is the unique
chord diagram with one chord.

Elements of the space $\cW_n=\cA^*_n$ (resp.\ $\ov{\cW}_n = \bcA_n^*$) 
are called {\em \ws s\/} (resp.{\em\ strong \ws s\/}) of order $n$.
We will view \ws s (resp. strong \ws s) of order $n$
as linear functions on $\langle \cD_n \rangle$
vanishing on the combinations of diagrams in 
 four-term relations~(\ref{eq:4term})  (resp.\ one-\ and
four-term relations). The spaces 
$$\cW=\oplus \cW_n \text{\quad and\quad } 
\ov{\cW} = \oplus \ov{\cW}_n
$$ 
are also $\Z$-graded Hopf algebras. The product of  $W_1 \in \cW_m$
and $W_2\in \cW_n$ is defined by
$$
(W_1\cdot W_2)(D) = \sum_{{E\subset D} \atop {|E|=m}}
		      W_1(E)W_2(D\setminus E),
\text{\quad where\ } D \in \cD_{m+n}.
$$ 
\bigskip

The fundamental result of the theory of \VI\ is the theorem
of Kontsevich \cite{Konts,BN} that the
Hopf algebra $\ov{\cW}$ is canonically isomorphic to the 
adjoint graded algebra of the filtered Hopf algebra
$\cV$. Let $\cV^f$ denote the space of \VI\ of framed knots.
Kontsevich proved that 
 $$
\cV_n/\cV_{n-1} \simeq \ov{\cW}_n,
    \text{\ \ and\ \ } 
\cV^f_n/\cV^f_{n-1} \simeq \cW_n~
$$
by explicitly constructing 
splitting maps
 \be
Z: \cW_n \to \cV^f_n \quad \text{and} \quad \ov{Z}: \ov{\cW}_n \to
\cV_n .
\label{eq:splitting}
\ee

The embedding $\bcA_n \hookrightarrow \cA_n$
gives a canonical {\em deframing projection\/}
\be
 \cW_n \ra \ov{\cW}_n
\label{eq:defr}
\ee
allowing to construct \VI\ of unframed knots from \ws\ that
give a priori only invariants of framed knots.

Many \ki s such as the Alexander and Jones polynomials are not \VI,
whereas their coefficients (after an appropriate change of variables)
are. Therefore, we can associate with such invariants not just one
\ws, but a sequence of \ws s $w_0, w_1, \ldots, w_n, \ldots$
where $w_n \in \cW_n$. 
We call elements of 
$$
\widehat{\cW} = \prod\limits_{n \ge 0} \cW_n
$$
{\em \Vas\ series\/} and write them as formal sums 
 $W=w_0+w_1+w_2 + \ldots$ (or sometimes as 
formal power series $W=\sum w_nz^n$, where $z$ is a formal parameter).
\Vas\ series can be viewed as linear functionals 
on the space generated by chord diagrams
$\cD=\bigoplus\limits_{n \ge 0} \langle \cD_n \rangle $.
For $D\in \cD_n$ and $W \in \cW$ we define $W(D) =w_n(D)$. 
The space $\widehat{\cW}$ also has an algebra structure
\be
(W\cdot W') (D) = \sum_{E\subset D} W(E)W'(D\setminus E).
\label{eq:product}
\ee

A \Vas\ series $W \in \widehat{\cW}$ is called {\em multiplicative\/} if
$$ 
 W(D_1\cdot D_2) = W(D_1)\cdot W(D_2) \ \text{\ for\ any\ }\ 
	 D_1, D_2 \in \cD.
$$

For our proof of the \MM\ conjecture it  is essential that, 
as we will see in Sections~4 and~7, both 
the colored Jones and the \AC\ \ws s are multiplicative.

\section{\Fd s and primitive elements of the Hopf algebra of diagrams} 

By the structure theorem on commutative and co-commutative Hopf
algebras~\cite{MiMo}, the algebra $\cA$ of diagrams is isomorphic to
the polynomial algebra
$S^\cdot (\cP)$, where $\cP$ is the space of primitive elements 
$$
\cP = \{ a\in \cA\ :\ \D(a) = 1\otimes a + a \otimes 1\}.
$$

The primitive space $\cP$ of $\cA$ can be very conveniently described
in terms of \Fd s.

A {\em \Fd \/} of order $n$ is a graph with $2n$ vertices of degrees 1
or 3 with a cyclic ordering on the set of its univalent ({\em
external\/}) vertices (or {\em legs\/})
and on each set of three edges meeting at a
trivalent ({\em internal\/}) 
vertex.
Drawing \Fd s we put legs on a circle (Wilson line) and
assume everything oriented in the counterclockwise direction.

\begin{rem}{\rm
\Fd s
are called  Chinese
Character diagrams in~\cite{BN}, but they are indeed \Fd s
arising in the perturbative Chern-Simons-Witten
quantum field theory~\cite{W,BN2}.
} \end{rem}

Let $\cF_n$ denote the set of all
\Fd s with $2n$ vertices (up to the natural equivalence of graphs with
orientations) having at least one external vertex on each connected
component. The set $\cD_n$ of chord diagrams with $p$ chords is a 
subset of $\cF_n$.
\medskip

\begin{prop}[\cite{BN}]
The embedding $\cD_n \hookrightarrow \cF_n$ induces an 
isomorphism of the space
$\cA_n$ of chord diagrams with $n$ chords modulo four-term
relations~(\ref{eq:4term}) and the space
$$
\langle\cF_n\rangle /\langle D_Y-D_{\mid\mid}+D_X\rangle
$$ 
of \Fn\ diagrams modulo the following
three-term relations
\begin{equation}\label{eq:3term}
\Picture{
\Arc[9]
\Arc[10]
\Arc[11]
\Arc[0]
\DottedArc[1]
\DottedArc[8]
\thinlines
\put(0.7,-0.7){\circle*{0.15}}
\qbezier(-0.1,0.1)(0.2,-0.5)(0.5,0.1)
\qbezier(0.2,-0.2)(0.3,-0.6)(0.7,-0.7)
\put(0,-1.4){\makebox(0,0){${}_{D_Y}$}}
}\quad = \quad
\Picture{
\Arc[9]
\Arc[10]
\Arc[11]
\Arc[0]
\DottedArc[8]
\DottedArc[1]
\Endpoint[10]
\Endpoint[11]
\thinlines
\put(0.5,-0.866){\line(-3,5){0.58}}
\put(0.866,-0.5){\line(-3,5){0.36}}
\put(0,-1.4){\makebox(0,0){${}_{D_{\mid\mid}}$}}
} \quad - \quad 
\Picture{
\Arc[9]
\Arc[10]
\Arc[11]
\Arc[0]
\DottedArc[8]
\DottedArc[1]
\Endpoint[10]
\Endpoint[11]
\thinlines
\put(0.5,0.1){\line(0,-1){0.97}}
\put(0.866,-0.5){\line(-5,3){1}}
\put(0,-1.4){\makebox(0,0){${}_{D_X}$}}
}~.
\end{equation}
\\[12pt]

In addition,
the following local relations hold for internal vertices of \Fd s
\begin{equation}
\begin{picture}(2,2)(0,0.375)
\qbezier(0.5,2)(1.65,1.3)(1.5,1)
\qbezier(1.5,2)(0.35,1.3)(0.5,1)
\qbezier(0.5,1)(1,0)(1.5,1)
\put(1,0){\line(0,1){0.5}}
\end{picture}
\quad = - 
\begin{picture}(2,2)(0,0.375)
\qbezier(0.5,2)(1,0)(1.5,2)
\put(1,0){\line(0,1){1}}
\end{picture}
\ {\mathrm and}\quad 
\begin{picture}(2,2)(-1,-1)
\put(0,-1.4){\line(1,1){1.4}}
\qbezier(-0.1,0.1)(0.2,-0.5)(0.5,0.1)
\qbezier(0.2,-0.2)(0.3,-0.6)(0.7,-0.7)
\end{picture} \quad =
\begin{picture}(2,2)(-1,-1)
\put(-0.034,-1.4){\line(1,1){1.4}}
\put(0.5,-0.866){\line(-3,5){0.58}}
\put(0.866,-0.5){\line(-3,5){0.36}}
\end{picture} \quad - 
\begin{picture}(2,2)(-1,-1)
\put(-0.034,-1.4){\line(1,1){1.4}}
\put(0.5,0.1){\line(0,-1){0.97}}
\put(0.866,-0.5){\line(-5,3){1}}
\end{picture}.\\[20pt]
\label{eq:asjac}
\end{equation}
\end{prop}
(Drawing relations between \Fd s we always assume that all
diagrams entering the same relation are identical except for the
depicted fragments.)
\medskip

This proposition implies, in particular, that
the space 
 $\langle\cF_n\rangle /\langle D_Y-D_{\mid\mid}+D_X\rangle$ 
is canonically dual to $\cW_n$ and that
any \ws\ $W:\cD_n \to k$ extends to $\cF_n$ by the 
rule 
$$
W(D_Y) = W(D_{||}) - W(D_X).
$$
\medskip

The realization of the algebra $\cA$ in terms  of \Fd s
endows it with a natural filtration consistent with the
grading
\be
\cA_n^{(0)} \subset \cA_n^{(1)} \subset \ldots \subset \cA_n^{(\ell)} \subset
\ldots  \cA_n^{(2n)} =\cA_n,
\ee
where $\cA_n^{(\ell)}$ is the subspace of $\cA_n$ generated by \Fd s
with $\le \ell$ legs. 
\bigskip

{}From the definition of the co-product in $\cA$ it follows
(see e.g.~\cite{BN}) that the primitive subspace $\cP\subset \cA$
is spanned by connected \Fd s and that the homogeneous component 
$\ov{\cP}_n$ of the primitive space of $\bcA$ coincides with
$\cP_n$ for $n\ne 1$, 
and $\ov{\cP}_1 = 0$ whereas $\cP_1=\langle \Th \rangle$.

Denote by $\cP_n^{(\ell)}$ the space $\cP \cap \cA_n^{(\ell)}$.
\bigskip

The following proposition leaves us with fewer diagrams in the
spanning set. Let us call a subgraph $T$ of a graph $F$ a {\em hanging
tree\/} if $T$ is a tree with more than one edge and becomes
disconnected from $F$ after removal of one of its leaves (terminal
vertices).
In particular, a hanging tree can be a connected component of $F$.
\begin{thm} \label{th:prim} \

(i) The space $\cP$ has a  basis consisting of connected \Fd s
without hanging trees.

(ii) For $n\ge 2$ the filtration 
\be
0 \subset
 \ldots \subset \cP_n^{(l-1)}\subset \cP_n^{(\ell)}	\subset \ldots
\subset  \cP_n
\label{eq:filtr}
\ee 
stabilizes at the  term $\cP_n^{(n)}$, i.e.
$$
\cP_n^{(n)} = \cP_n^{(n+1)} = \ldots = \cP_n.
$$

(iii) The leading quotient space 
$\cP_n^{(n)}/\cP_n^{(n-1)}$ of~(\ref{eq:filtr}), $n \ge 2,$
is generated by the $n$-spoke wheel diagram
\be
w_n = 
\Picture{
\DottedCircle
\put(0,0){\circle{0.6}}
\put(1,0){\line(-1,0){0.7}}
\put(0.707,0.707){\line(-1,-1){0.5}}
\put(0.707,-0.707){\line(-1,1){0.5}}
\put(0,1){\line(0,-1){0.7}}
\put(0,-1){\line(0,1){0.7}}
\bezier{10}(-0.5,0.5)(-0.9,0)(-0.5,-0.5)
\Arc[3]\Arc[2]\Arc[1]\Arc[0]\Arc[11]\Arc[10]\Arc[9]
}
\ee\\[10pt]
\end{thm}

First we will prove several simple lemmas.

\begin{lemma} \label{lem:bubble}
The following identities hold in the algebra $\cA$ of
diagrams 
\begin{equation}
\Picture{
\Arc[8]
\Arc[9]
\Arc[10]
\DottedArc[7]
\DottedArc[11]
\thinlines
\put(0,-0.4){\circle{0.4}}
\put(0,-1){\line(0,1){0.4}}
\put(0,-0.2){\line(0,1){0.5}}
} 
= 2~\Picture{
\Arc[8]
\Arc[9]
\Arc[10]
\DottedArc[7]
\DottedArc[11]
\thinlines
\qbezier(-0.29,-0.96)(0,-0.25)(0.29,-0.96)
\put(0,-0.6){\line(0,1){0.9}}
} 
\end{equation}\\[10pt]
and\\[8pt]
\begin{equation}
\Picture{
\put(0,0){\line(1,-1){0.7}}
\put(-1,0){\line(1,0){0.3}}
\put(0,0){\line(0,1){1}}
\put(-0.5,0){\circle{0.4}}
\put(0,0){\line(-1,0){0.3}}
} 
= 2\ \Picture{
\put(0,0){\circle{0.4}}
\put(0,1){\line(0,-1){0.8}}
\put(-0.2,0){\line(-1,0){0.8}}
\put(0.141,-0.141){\line(1,-1){0.5}}
} 
=
\Picture{
\put(0,1){\line(0,-1){0.3}}
\put(0,0){\line(1,-1){0.7}}
\put(0,0){\line(-1,0){1}}
\put(0,0.5){\circle{0.4}}
\put(0,0){\line(0,1){0.3}}
}~. \label{eq:w3}
\end{equation}
\end{lemma}
\begproof
These identities follow immediately from the three-term
relations~(\ref{eq:3term}) and from~(\ref{eq:asjac}).  
\eproof

\begin{lemma} \label{prop:trees}
Let $F$ be a \Fd\ with $\ell$ legs and a planar
hanging tree $T$ with $d$
legs that does not intersect the edges of $F\setminus T$.
Then 
$$
F= 2^{-(d-1)} F',
$$
where $F'$ is the diagram with $\ell- d + 1$ legs obtained
from $F$ by replacing the hanging tree $T$ with the chain 
of $d-1$ ``bubbles''. For example,  
$$ 
\Picture{
\put(0.2,0){\line(-1,0){0.3}}
\put(0.2,0){\line(1,1){0.6}}
\put(0.2,0){\line(1,-1){0.6}}
\put(0.453,-0.253){\line(1,1){0.5}}
\DottedArc[2]\DottedArc[10]
\Arc[1]\Arc[0]\Arc[11]
} \
= \frac{1}{4}\cdot
\Picture{
\put(1,0){\line(-1,0){0.3}}
\put(0.3,0){\line(-1,0){0.3}}
\put(-0.4,0){\line(-1,0){0.3}}
\put(0.5,0){\circle{0.4}}
\put(-0.2,0){\circle{0.4}}
\DottedArc[2]\DottedArc[10]
\Arc[1]\Arc[0]\Arc[11]
}
$$
\end{lemma}

\begproof
The statement follows from Lemma~\ref{lem:bubble} by induction on $d$.
\eproof

\begin{lemma} \label{lem:tree}
A \Fd\ $F$ with $\ell$ legs and a hanging tree belongs to 
$\cA^{(\ell-1)}$. 
\end{lemma}

\begproof
Repeatedly applying the three-term relation~(\ref{eq:3term})
and the relation~(\ref{eq:asjac})
we can transform $F$ to a \Fd\ $\widetilde{F}$ with a planar hanging tree 
not intersecting the rest of $\widetilde{F}$ so that
$F = \pm\widetilde{F} + G$, where  $G$ is a linear combination
of diagrams with $\le \ell -1$ legs. Therefore $G\in \cA^{(\ell -1)}$
and, by~\ref{prop:trees}, $\widetilde{F}\in \cA^{(\ell -1)}$.
\eproof

\begin{lemma} \label{lem:betti}
Let $F$ be a connected \Fd\ of order $n$ with
$\ell$ legs. Then the first Betti number of $F$ is
$$
h^1(F) = n - \ell +1.
$$
\end{lemma}
\begproof
The diagram $F$ has $2n-\ell$ trivalent vertices and, consequently,
$e = \half(\ell +3(2n-\ell)) = 3n-\ell$ edges. Since $F$ is connected,
its Euler characteristic is
$\chi(F)=1 - h^1(F) = 2n - e = \ell - n$, and 
$h^1(F)=n-\ell +1$.
\eproof
\medskip

\noindent
{\em Proof of Theorem~\ref{th:prim}.}
Choose a basis $\cB^{(\ell)}$ 
in $\cP_n^{(\ell)}$ represented by connected \Fd s
with $\le \ell$ legs such that $\cB^{(\ell)} \subset \cB^{(\ell+1)}$.
If $F \in \cB^{(\ell +1)} \setminus \cB^{(\ell)}$ then
it has exactly $\ell$ legs and, by Lemma~\ref{lem:tree}, cannot have
hanging trees. Therefore, 
$\bigcup\limits_{\ell \ge 0} \cB^{(\ell)}$ is a 
basis of $\cP_n$ satisfying {\em (i).}

Let $m$ be the largest $\ell$ for which
$\cP_n^{(\ell)} \ne \cP_n^{(\ell-1)}$ and
$F\in \cB_n^{(m)} \setminus \cB_n^{(m-1)}$.
Since $F$ is a connected  
diagram with $m$ legs
and without hanging trees and $n \ge 2$, the 
diagram  $F$ cannot be a tree and its first Betti number 
$h^1(F)$ must be at least $1$. By Lemma~\ref{lem:betti}
$h^1(F)=n-m +1 \ge 1$ from where $m \le n$. This gives
part {\em (ii)\/} of the theorem.

In the case $F \in \cB^{(n)}$ the diagram $F$ has exactly one cycle
and, since $F$ does not contain hanging trees, all the $n$ legs of $F$ 
are adjacent to this cycle. Therefore, $F$ is an $n$-spoked wheel
diagram with
a possibly different ordering of legs than that of the standard wheel
diagram $w_n$.
By applying relations~(\ref{eq:3term}) and~(\ref{eq:asjac})
we can
change the ordering of legs of $F$ at the cost of adding extra
diagrams with {\em less\/} than $n$ legs.
Therefore 
% $$
% F = \pm w_n + \sum_i F_i, \text{\ \ where \ } F_i \in \cP_n^{(n-1)}
% $$ and
 $F = \pm w_n$ modulo $\cP_n^{(n-1)}$ which proves {\em (iii).}
\eproof

\begin{rem}{\rm It follows from (iii) that 
\ $\text{dim} \cP_n^{(n)}/\cP_n^{(n-1)} \le 1.$
Chmutov and Varchenko proved in~\cite{CV} a stronger statement
that this dimension is equal to one if and only if $n$ is even,
but we will not need this for the proof of the \MM\ conjecture. 
}\end{rem}

\section{Weight systems coming from Lie algebras} 
\label{sec:FG}

Here we recall a construction that
assigns a \Vas\ series to every Lie
(super)\-algebra with an invariant inner product.
(For more details see~\cite{VaiYB}.)
\bigskip

Let $L$ be a Lie (super)algebra with an $L$-invariant inner product
$b: L\otimes L \ra k$. To each \Fn\ diagram $F$ with $m$ univalent
vertices  we assign a tensor 
$$
   T_{L,b}(F) \in L^{\otimes m}
$$ 
as follows.

The Lie bracket $[\ ,\ ]: L\otimes L \to L$ can be considered
as a tensor in $L^*\otimes L^* \otimes L$. The inner product 
$b$ allows us to identify the $L$-modules 
$L$ and $L^*$, and therefore  $[\ ,\ ]$ can be considered as a tensor
$f\in (L^*)^{\otimes 3}$ and $b$ gives rise to an invariant symmetric
tensor $c \in V\otimes V$.

For a \Fn\ diagram $F$ denote by $T$ the set of its trivalent
vertices, by $U$ the set of its univalent (exterior) vertices, and by
$E$ the set of its edges. Taking $|T|$ copies of the tensor $f$ and $|E|$
copies of the tensor $c$ we construct a new tensor
$$
\widetilde{T}_L(F) = \Bigl(\bigotimes_{v\in T} f_v\Bigr) \otimes
		     \Bigl(\bigotimes_{\ell\in E} c_\ell\Bigr)
$$
which is an element of the tensor product
$$
{\cal{L}}^F = 
\Bigl( \bigotimes_{v \in T}
(L^*_{v,1}\otimes L^*_{v,2} \otimes L^*_{v,3})
\Bigr)
\otimes
\Bigl( \bigotimes_{\ell \in E}
(L_{\ell,1}\otimes L_{\ell,2})
\Bigr),
$$
where $(v,i), \ i=1,2,3,$ mark the three edges meeting at the
vertex $v$ (consistently with the cyclic ordering of these edges), and
$(\ell,j), \ j=1,2,$ denote the endpoints of the edge $\ell$. 

Since $c$ is symmetric and $f$ is completely antisymmetric,
the tensor $\widetilde{T}_L(F)$ does not depend on the choices of
orderings.

If $(v,i)=\ell$ and $(\ell,j)=v$, there is a natural contraction map
$L^*_{v,i} \otimes L_{\ell,j} \to k$. Composition of all such 
contractions gives us a map
$$
{\cal{L}}^F \longrightarrow  \bigotimes_{u \in U} L = L^{\otimes m}, \ \text{where} \ m=|U|.
$$

The image of $\widetilde{T}_L(F)$ in
$L^{\otimes m}$ is denoted by $T_{L,b}(F)$ (or usually just by
$T_{L}(F)$).

(In the case of Lie superalgebras we also have to take special care of
signs.  See \cite{VaiYB} for details.)

For example, for the diagrams 
\def\DiagC{
\begin{picture}(3,1)
\qbezier(0.5,0)(1.5,2)(2.5,0)
\thicklines
\put(0,0){\line(1,0){3}}
\end{picture}
}
\def\Bubble{
\begin{picture}(4,1.5)
\put(1.5,0){\oval(2,2)[tl]}
\put(2.5,0){\oval(2,2)[tr]}
\put(2,1){\circle{1}}
\thicklines
\put(0,0){\line(1,0){4}}
\end{picture}}
\def\DiagramK{
\begin{picture}(5,1.5)
\qbezier(0.5,0)(1.25,1.5)(2,0)
\qbezier(3,0)(3.75,1.5)(4.5,0)
\qbezier(1.25,0.75)(2.5,2)(3.75,0.75)
\thicklines
\put(0,0){\line(1,0){5}}
\end{picture}}
\be
C= \DiagC~, 
\quad 
B = \Bubble~, 
\quad \text{and} 
\quad
K =\DiagramK
\label{eq:DiagK}
\ee
we have 
$$T_L(C) = \sum_{ij}  b^{ij}e_i\otimes e_j= c,
$$ 
the Casimir element corresponding to the inner product~$b$, \
$$
T_L(B) = \sum b^{is}b^{tj}b^{kp}b^{lq}f_{skl}f_{pqt}e_i \otimes e_j,
$$
the tensor in $L\otimes L$
corresponding to the Killing form on $L$ under the identification 
$L^* \simeq L$, and
$$
T_L(K) = \sum b^{in} b^{jp} b^{qr} b^{kt} b^{\ell s}
f_{npq} f_{tsr} e_i \otimes e_j
\otimes e_k \otimes e_\ell~,
$$
 where $f^i_{jk}=\sum b^{is}f_{sjk}$ are the structure constants of 
$L$ in a basis $e_1,e_2,\ldots$~.

Tensor $T_L(F)$ is invariant with respect to the $L$-action on
$L^{\otimes m}$ and its image $V_L(F)$ in the universal enveloping
algebra $\cU(L)$ belongs to the center 
$Z(\cU(L))=\cU(L)^L$ and does not depend on the place where we
cut the Wilson line to obtain a linear ordering of the external
vertices of $F$.

Therefore, for every Lie algebra $L$ with an invariant inner product
there exists a natural \Vas\ series
$$
V_L: \cD \to Z(\cU(L)) 
$$
which is called {\em the universal weight system}
corresponding to $L$ and $b$.
It is universal in the sense that any \Vas\  series $V_{L,R}$
constructed using a representation $R$ of the Lie algebra $L$ (see
\cite{BN})  is an evaluation of $V_L$:
$$
V_{L,R}(D) = \text{Tr}_{R}\bigl(V_L(D)\bigr).
$$

\begin{prop}
Both the universal
 \Vas\ series $V_L$ and its deframing $\ov{V}_L$ are
 multiplicative, i.e. 
$$V_L(D_1 \cdot D_2) = V_L(D_1)V_L(D_2)
\text{\ \ and \ \ }
\ov{V}_L(D_1 \cdot D_2) = \ov{V}_L(D_1)\ov{V}_L(D_2) .
$$
\end{prop}
\begproof
Multiplicativeness of $V_L$ follows immediately from its construction. 
Deframing $\ov{W}$ of a multiplicative \Vas\ series $W$ is 
again multiplicative, since the embedding $\bcA_n \hookrightarrow
\cA_n$ is a homomorphism of algebras and
$W$ and $\ov{W}$ coincide on all
primitive \Fd s $F$ in $\cP$ but one
$$
\ov{W}(F) = \cases{W(F) & if $F \ne \Th$, 
                                           \cr 
		    0 & if $F=\Th$.}
$$
Therefore $\ov{V}_L$ is also multiplicative.
\eproof

\begin{corol}
The colored Jones and the \AC\ \ws\ are multiplicative.
\end{corol}
\begproof
This is so because the deframed
universal  $sl_2$ and \gl\ 
\ws s are the \Vas\ series corresponding to
the colored Jones and \AC\ invariants respectively
(see Sections~5 and~7).
\eproof

\section{Reduction of the \MM\ conjecture to \Fd s}
\label{sec:canon}

There are two constructions of knot invariants from 
a semi-simple Lie algebra $L$  with an invariant inner 
product $b$ and a representation $R$. 
The Reshetikhin-Turaev construction~\cite{RT} 
based on quantum groups gives
invariants $I_{L,R}$ and $\bar{I}_{L,R}$
(of framed and unframed knots, resp.)\ with values
in $\Z[t,t^{-1}]$. The second way is to apply Kontsevich's
construction to the
\Vas\ series $V_{L,R}=\text{Tr}_R(V_L)$ and $\ov{V}_{L,R}$ discussed  
in the previous section. 

These two 
constructions are equivalent. 
The coefficients of the power series expansions of 
$\bar{J}_{L,R}(z)=\bar{I}_{L,R}(e^z)$ are \VI,
 and the corresponding series of
\ws s coincide with $\ov{V}_{L,R}$ (see~\cite{BL,P}).
Conversely, Kontsevich's splitting maps~(\ref{eq:splitting})
applied to $\ov{V}_{L,R}$ and $V_{L,R}$)
gives the sequences of the coefficients of $\bar{J}_{L,R}$
and $J_{L,R}$ (see~\cite{Kassel,LM}).

A \Vas\ invariant is called {\em canonical\/} 
if it belongs to the image of 
Kontsevich's map $Z: \cW_n \to \cV_n$~.
A formal power series  
$\sum\limits_{n \ge 0} a_n z^n \in \cV[[z]]$ 
is called {\em canonical\/} if every coefficient $a_n$ 
is a canonical \Vas\ invariant of order $\le n$.

The colored Jones invariant is the canonical invariant
$Z(\frac{1}{d}\ov{W}_{sl_2,R_d})$ corresponding to the $d$-dimensional
representation of $sl_2$.

\BN\ and \Gar~\cite{BNG} proved that
\be
\widetilde{\D}(z) =  \frac{z}{e^{z/2} - e^{-z/2}}\D(e^z),
\label{eq:alexnorm}
\ee
where $\D(K,t)$ is the Alexander polynomial of knot $K$,
is a canonical series whose \Vas\ series
coincides with the series $C$ of the \AC\ polynomial.

A canonical invariant is uniquely determined by its \ws\ and
the product of two canonical invariants or \Vas\ series
$Z(W_1)$ and $Z(W_2)$ is again canonical with the \ws\
equal to $W_1\cdot W_2$ (see~\cite{BNG}).
Therefore, to prove the \MM\ conjecture 
it is enough to establish the relations~(\ref{eq:MM0})
and~(\ref{eq:MM}) on the level of \ws s.

Theorem~\ref{thm:mm} is thus reduced to the following
relations between \ws s.

\begin{prop} \label{prop:ws}
Let 
\be
{\hW}_{sl}=\frac{1}{d}\sum_{n \ge 0} \ov{W}_{sl_2,R_d,n}
=   \sum_{i,n \ge 0} w_{in} d^i
\label{eq:coefmm}
\ee
be  the deframed \Vas\ series coming
from the $d$-dimensional representation of $sl_2$ with the standard
metric $\langle x,y \rangle = \text{Tr}(xy)$
normalized by dividing by $d$.
Then
\be (i) \qquad w_{in}=0 \  \text{~for~} \  i>n
\label{eq:mm1}
\ee
and
\be (ii) \qquad (\sum_{n \ge 0} w_{nn} )\cdot C = \e,
\label{eq:mm2}
\ee
where $C = \sum C_n z^n \in \widehat{\cW}$ is the \Vas\ series
of the \AC\ invariant and $\e$ is the \Vas\ series 
$$
\e(D) = \cases {1 & if $|D|=0$ ,\cr 0 & if $|D| > 0$~. }
$$
\end{prop}
\medskip

We will derive~(\ref{eq:mm1}) and~(\ref{eq:mm2}) 
from Theorem~\ref{th:prim} and the following two propositions.

Let $V_{sl}$ be the universal \Vas\ \ws\ corresponding to the Lie algebra
$sl_2$ with the standard invariant form $\langle x,y\rangle =
\text{Tr}(xy)$. We consider $V_{sl}$ as a homomorphism from
the algebra of diagrams to
 $\C[c] = Z(\cU(sl_2))$, where $c$ is the quadratic Casimir.
Since $V_{L}$ respects the filtrations in $\cA$ and $\cU(L)$,
we have 
\be \label{eq:deg}
\text{deg}_c V_{sl}(F) \le \left[\frac{\ell}{2}\right] \text{\ \ for \
\ } F \in \cP^{(\ell)}.
\ee
\begin{prop} \label{prop:slprim}
Let $w_{2n}$ be the wheel diagram with $2n$ spokes, $n\ge 2$.
Then 
$$
V_{sl}(w_{2n}) = 2^{n+1} c^n + \text{\ lower\ order\ terms\  in\ } c
$$
\end{prop}

This  proposition will be proved in Section~6.

\begin{prop} \label{prop:glprim}
Let $C$ be the \AC\ \ws. Then for a connected \Fd\ $F$
of order $n\ge 2$ we have 
\be \label{eq:ac}
C(F) = \cases{0 & if $F \in \cP_n^{(n-1)}$,\cr
             -2 & if $F=w_n$ and $n$ is even,\cr 
              0 & if $F=w_n$ and $n$ is odd.}
\ee
\end{prop}

This proposition is proved in Section~7.
\bigskip

\noindent
{\em Proof of the \MM\ conjecture (Proposition~\ref{prop:ws}).\/}
Let 
$$ 
\ov{V}_{sl}(c) = \sum_{n\ge 0} \ov{V}_n(c) = \sum_{n,i \ge 0}
\ov{v}_{in} c^i, 
$$
be 
the deframed universal $sl_2$ \Vas\ series considered as a function on
$c$ and $\ov{V}_n$ be its order $n$ component.

The quadratic Casimir $c\in Z(\cU(sl_2))$ in the $d$-dimensional
irreducible representation acts as a multiplication by $(d^2-1)/2$.
Therefore, 
\be
\frac{1}{d} \ov{W}_n(d) = \ov{V}_n\left(\frac{d^2-1}{2}\right),
\label{eq:susbt}
\ee
and equations~(\ref{eq:mm1}) and~(\ref{eq:mm2})  
become equivalent to
\bea
 (i)& & \qquad \ov{v}_{in}=0 \  \text{~for~} \  i>n/2 \qquad \
\label{eq:mmu1} \\
&\text{and} & \nonumber\\
(ii)& & \qquad \widetilde{V}\cdot C = \e, \qquad \
\label{eq:mmu2}
\eea
where 
$
  \widetilde{V}=\sum\limits_{n \ge 0} 2^{-n}\ov{v}_{n,2n}.
$

Since both ${\ov{V}}_{sl}$ and $C$ are multiplicative \Vas\ series,
it is enough to verify~(\ref{eq:mmu1})
and~(\ref{eq:mmu2}) on primitive elements of $\cA$.
On diagrams of degree $0$ or $1$ these equations are trivially
satisfied.

For $n \ge 2$ by Theorem~\ref{th:prim}.(ii) 
the primitive space 
$\cP_n$ is spanned by \Fd s with  $\le n$ legs. Therefore 
by~(\ref{eq:deg})
$\text{deg}_c V_{sl}(F) \le n/2$
for any $F \in \cP_n$ and
$\ov{v}_{in}(F)=0$ for $i > n/2$. This gives~(\ref{eq:mmu1})
and the first part of the \MM\ conjecture.

Since $V_{sl}$ and $\ov{V}_{sl}$ coincide on
$\cP_n$ for $n\ge 2$, we have by~(\ref{eq:deg}) and
Proposition~\ref{prop:slprim}
\be \label{eq:vtilde}
\widetilde{V}(F) =  
\cases{0 & if $F \in \cP_n^{(n-1)}$,\cr
              2 & if $F=w_n$ and $n$ is even,\cr 
              0 & if $F=w_n$ and $n$ is odd.}
\ee

If $F\in \cP_n$, then by the definition of a primitive element
\begin{eqnarray}
(\widetilde{V} \cdot C)(F) &=& (\widetilde{V} \otimes C) \D(F)
\nonumber\\ 
& = & (\widetilde{V} \otimes C) (F\otimes 1 + 1 \otimes F) \nonumber\\ 
& = & \widetilde{V}(F) + C(F). \label{eq:vanish}
\end{eqnarray}

If $F\in \cP_{n}^{(n-1)}$
then both terms in~(\ref{eq:vanish}) vanish due 
to~(\ref{eq:ac}) and~(\ref{eq:vtilde}). 
Therefore, it only
remains to check that $\widetilde{V}(F) + C(F) = 0$ for
$F\in  \cP_{n}^{(n)} \setminus \cP_{n}^{(n-1)}$. 
In this case according
to Theorem~\ref{th:prim}.(iii) we can assume that $F=w_n$ and
again by~(\ref{eq:ac}) and~(\ref{eq:vtilde})
we get $\widetilde{V}(F)= -C(F)$.
\eproof

\section{Proof of Proposition~5.2}

\def\DiagC0{
\begin{picture}(3,1)
\qbezier(0.5,0)(1.5,2)(2.5,0)
\end{picture}
}
\def\bubble{
\begin{picture}(4,1.5)
\put(1.5,0){\oval(2,2)[tl]}
\put(2.5,0){\oval(2,2)[tr]}
\put(2,1){\circle{1}}
\end{picture}}
\def\DiagramK0{
\begin{picture}(4.5,1.5)
\qbezier(0.5,0)(1.25,1.5)(2,0)
\qbezier(3,0)(3.75,1.5)(4.5,0)
\qbezier(1.25,0.75)(2.5,2)(3.75,0.75)
\end{picture}}

\begin{prop}
The universal $sl_2$ \ws\ $V_{sl}$ satisfies the
relations

\bea
(i) &&  V_{sl}\Biggl(\bubble \Biggr) = 
4 V_{sl}\Biggl(\DiagC0 \Biggr)     \label{eq:slkilling}\\[12pt]
& \text{and} &\nonumber\\[10pt]
(ii) && V_{sl}\Biggl(\DiagramK0 \Biggr) = 2V_{sl} \Biggl(
~\begin{picture}(4,1.5)
\qbezier(0.5,0)(2,2.5)(3.5,0)
\qbezier(1.5,0)(2,1)(2.5,0)
\end{picture}
-~
\begin{picture}(4,1.5)
\qbezier(0.5,0)(1.5,2)(2.5,0)
\qbezier(1.5,0)(2.5,2)(3.5,0)
\end{picture}
\Biggr)  \label{eq:lagrange}
\eea
\end{prop}

\begproof
The standard invariant metric on $L=sl_2$ is 
$\langle x,y \rangle = \text{Tr}_V (xy). $
This metric is equal to the one fourth of the  Killing form which is the
tensor represented by the `bubble' diagram in~(\ref{eq:slkilling}). 
This gives (i).

The relation~(\ref{eq:lagrange}) is a relation between invariant
tensors in $L^{\otimes 4}$ which is equivalent to the following
classical identity for the bracket and the inner product of $L$
\be 
\langle[a,b],[c,d]\rangle = 2\langle a,d \rangle \cdot \langle b,c
\rangle  -2 \langle a,c \rangle \langle b,d \rangle
\ \text{ \ for \ } \ a,b,c,d \in sl_2 .
\label{eq:alglagr}
\ee

This relation is better known in the $so_3$ realization of the Lie
algebra $L=sl_2$.
Under the isomorphism  $sl_2 \simeq so_3$, the inner product and the Lie
bracket on $sl_2$ become respectively the scalar and the vector product in
the three-dimensional space and~(\ref{eq:alglagr}) becomes the classical
Lagrange's identity%
\footnote{In fact, the identity~(\ref{eq:class}) goes back to
Euler's theory of 
the motion of a rigid body. In the notation of 20th century physics it
looks like 
{$\epsilon_{\a\b\g}\epsilon_{\rho\g\s} =\d_{\a\s}\d_{\b\rho} -
\d_{\a\rho}\d_{\b\s}, 
$}
where $\d$ is the Kronecker delta, and $\e$ is the standard completely
antisymmetric tensor in $\R^3$.}{}
\be \label{eq:class}
\mathbf 
[a\times b]\cdot [c\times d] = (a\cdot c)(b\cdot d) - (a\cdot
d)(b\cdot c)
\ee
which is equivalent to the better known {\em fundamental relation of
vector calculus\/} 
$$\mathbf
[a \times [b\times c]] = (a\cdot c) b - (a\cdot b) c~.
$$
\eproof

{\begin{rem}{\rm
Identity~(\ref{eq:lagrange}) allows us to compute the
universal $sl_2$ framed and unframed \ws s recursively.
See~\cite{CV} and~\cite{Vmm}  for the recursion formulas.
} \end{rem}}

\noindent
{\em Proof of Proposition~\ref{prop:slprim}.\/}
Applying~(\ref{eq:lagrange}) 
to the wheel diagram $w_{2n}$ and using
the three-term relation we get
$$
V_{sl}(w_{2n}) = V_{sl}\Bigl(2\Th \cdot w_{2n-2} - 2 T + 2w_{2n-1}\Bigr),
$$

% $$ 
%    \fbox{FIGURES ???}
% $$
where $T$ is a planar tree with $2n$ legs. According to
Lemma~\ref{prop:trees},  $T$ is proportional to a diagram with
only $2$ legs, therefore both $V_{sl}(T)$ and $V_{sl}(w_{2n-1})$
are polynomials in $c$ of degree less than $n$ and do not
contribute to the degree $n$ term in $V_{sl}(w_{2n})$.
Therefore, keeping track of only highest degree terms we get
$$
V_{sl}(w_{2n}) = 2c V_{sl}(w_{2n-2}) + \ldots =
(2c)^{n-1}V_{sl}(w_2) = 2^{n+1} c^n,
$$
since by~(\ref{eq:slkilling})
 $V_{sl}(w_2)= 4c$.
\eproof

\section{\AC\  \ws}
Here we prove Proposition~\ref{prop:glprim}. It follows from special
relations satisfied by the \AC\ \ws\ $C$.

\begin{prop}\label{prop:gl}
The \AC\ \Vas\ series $C$ satisfies the following relations
\bea
(i)&& 
C\Biggl( 
%  \fbox{Bubble w. two bold ends} 
\Picture{
\DottedCircle
\put(-1,0){\line(1,0){0.8}}
\put(1,0){\line(-1,0){0.8}}
\put(0,0){\circle{0.4}}
\Arc[6]\Arc[0]
}
\Biggr)
= -2 C\Biggl(
% \fbox{Bubble removed}
\Picture{
\DottedCircle
\Arc[6]\Arc[0]
}
\Biggr),
\label{eq:ACi}
\\[20pt]
(ii)&& 
C\Biggl(
%  \fbox{Bubble connected to a 3-valent vertex}
\Picture{
\put(0,1){\line(0,-1){0.3}}
\put(0,0){\line(1,-1){0.7}}
\put(0,0){\line(-1,-1){0.7}}
\put(0,0.5){\circle{0.4}}
\put(0,0){\line(0,1){0.3}}
} 
\Biggr) = 0,\label{eq:ACii}
\\[20pt]
(iii) &&
C\Biggl(
\Picture{
\put(0,0){\line(0,1){0.56}}
\put(0,0){\line(1,-1){0.5}}
\put(0,0){\line(-1,-1){0.5}}
\put(-0.5,-0.5){\line(0,-1){0.5}}
\put(-0.5,-0.5){\line(-1,0){0.5}}
\put(0.5,-0.5){\line(0,-1){0.5}}
\put(0.5,-0.5){\line(1,0){0.5}}
\put(0,0.56){\line(1,1){0.5}}
\put(0,0.56){\line(-1,1){0.5}}
} 
\Biggr) = 0,\label{eq:ACiii}
\\[20pt]
(iv)&& 
C\Biggl(\DiagramK0 \ \Biggr) \nonumber \\
& & = 
\half C\Biggl(%
\begin{picture}(3,1)
\put(0.75,0){\oval(1,1)[tl]}
\put(1.25,0){\oval(1,1)[tr]}
\put(1,0.5){\circle{0.5}}
\put(2,0){\oval(1.5,1)[t]}
\end{picture} + 
\begin{picture}(3,1)
\put(1,0){\oval(1.5,1)[t]}
\put(1.75,0){\oval(1,1)[tl]}
\put(2.25,0){\oval(1,1)[tr]}
\put(2,0.5){\circle{0.5}}
\end{picture} - 
\begin{picture}(3,1.5)
\put(1.25,0){\oval(1,1)[tl]}
\put(1.75,0){\oval(1,1)[tr]}
\put(1.5,0.5){\circle{0.5}}
\put(1.5,0){\oval(2.5,2)[t]}
\end{picture} -
\begin{picture}(2.5,1.5)
%\begin{picture}(3,1.5)
\put(1.5,0){\oval(1.5,1)[t]}
\put(1.25,0){\oval(2,2)[tl]}
\put(1.75,0){\oval(2,2)[tr]}
\put(1.5,1){\circle{0.5}}
\end{picture}\
\Biggr).\label{eq:ACiv}
\eea
\end{prop}
\bigskip

\noindent
{\em Proof of Proposition~\ref{prop:glprim}.\/}
Let $F \in \cP_n^{(n-1)}$ be a connected \Fd\ with 
$\ell \le n-1$ legs. Since the number of
trivalent vertices in $F$ adjacent to the Wilson line
cannot exceed the number of legs and $\ell < n <2n - \ell$,
there is at least one trivalent vertex in $F$ adjacent
to only trivalent vertices. Therefore $F$ contains a fragment
shown on~(\ref{eq:ACii}) or~(\ref{eq:ACiii}), and by
Proposition~\ref{prop:gl}.(ii) and (iii), $C(F)=0$.

To evaluate $C$ on the wheels $w_n$ 
we apply (iv) to $w_n$ and taking into account (i),  (ii)
and~(\ref{eq:w3}) we obtain
$$
C(w_n) = C(w_{n-2}) = \ldots = 
\cases{ C(w_2) = - 2 & if $n$ is even,\cr
	C(w_3) = 0 & if $n$ is odd.}
$$
% 
% $$ 
%    \fbox{FIGURES ???}
% $$
\eproof
\bigskip

\noindent
{\em Proof of Proposition~\ref{prop:gl}.\/}
The
Lie superalgebra $L=\gl$
of endomorphisms of the $(1|1)$-dimensional
superspace has an invariant inner product
$\langle x,y\rangle = str (xy).$ 

Therefore, we can consider the universal \ws\
$V_{gl(1|1)}$ with values in 
$$Z(\cU(\gl)) = \C[h,c],
{\text \ \ where\ }
h = \pmatrix{1 & 0\cr 0 & 1\cr} \in \gl 
\text{ \ and \  } c \text{ \ is\ the\ quadratic\ Casimir.}
$$

As it was proved in~\cite{FKV}, the deframed \Vas\ series
$\ov{V}_{L}$ coincides with the \AC\ \ws.
Namely, for $F \in \cF_n$
$$ \ov{V}_L(F) = C(F)h^{n/2},$$
in particular, $\ov{V}$ vanishes on $\cF_n$ for odd $n$.

In~\cite{FKV} we proved
the following analogs of the fundamental
$sl_2$ relations~(\ref{eq:slkilling}) and~(\ref{eq:lagrange})
between invariant tensors on $L=\gl$ defined by the diagrams
$B$ and $K$ (see (\ref{eq:DiagK})) 

\be
V_L \Biggl( 
\Picture{
\put(-1,0){\line(1,0){0.8}}
\put(1,0){\line(-1,0){0.8}}
\put(0,0){\circle{0.4}}
}
\Biggr)
=  - 2 h^2
\label{eq:gl1bb}
\ee
and
\be
V_L\Biggl( \DiagramK0 \Biggr) 
= \half V_L(M),
\label{eq:gl11fund}
\ee
where
$$
M =
\begin{picture}(3,1)
\put(0.75,0){\oval(1,1)[tl]}
\put(1.25,0){\oval(1,1)[tr]}
\put(1,0.5){\circle{0.5}}
\put(2,0){\oval(1.5,1)[t]}
\end{picture} + 
\begin{picture}(3,1)
\put(1,0){\oval(1.5,1)[t]}
\put(1.75,0){\oval(1,1)[tl]}
\put(2.25,0){\oval(1,1)[tr]}
\put(2,0.5){\circle{0.5}}
\end{picture} - 
\begin{picture}(3,1.5)
\put(1.25,0){\oval(1,1)[tl]}
\put(1.75,0){\oval(1,1)[tr]}
\put(1.5,0.5){\circle{0.5}}
\put(1.5,0){\oval(2.5,2)[t]}
\end{picture} -
\begin{picture}(3,1.5)
\put(1.5,0){\oval(1.5,1)[t]}
\put(1.25,0){\oval(2,2)[tl]}
\put(1.75,0){\oval(2,2)[tr]}
\put(1.5,1){\circle{0.5}}
\end{picture}~.
$$
\medskip

These identities give relations~(\ref{eq:ACi}) and~(\ref{eq:ACiv})
for the \AC\ \ws.

The identity~(\ref{eq:gl1bb}) implies~(\ref{eq:ACii}),
since $h$ is a central element in \gl.
Finally, the relation~(\ref{eq:ACiii})  follows from~(\ref{eq:ACii})
and~(\ref{eq:ACiv}). Algebraically it is equivalent to the well-known
fact that \gl\ is solvable of depth three, i.e. that 
\be
\Bigl[\bigl[ [a,b],[c,d]\bigr],g\Bigr] = 0, \text{\ \ for \ any \ \ }
a,b,c,d, g \in \gl .
\label{eq:gl1solv}
\ee
Equation~(\ref{eq:gl1solv})
is obvious from a glance at the multiplication table for
the bracket on \gl\ which shows that
 $\bigl[ [a,b],[c,d]\bigr]$ is always a multiple of the central element
$h$. 
\eproof

\begin{rem}
{\rm
A direct combinatorial proof of the key Proposition~\ref{prop:gl}
which is not based on the relation between the \AC\ 
polynomial and Lie superalgebra \gl\
follows from our lemmas in Section~4 of \cite{FKV}. A different proof
of Proposition~\ref{prop:glprim} can be derived from the
characterization of the space generated by the \AC\ \ws s\
in~\cite{AK}. 
}\end{rem}

\bigskip

\noindent
{\bf Acknowledgments.}
S.~Chmutov informed me that he has also found a proof of the \MM\
conjecture based on similar ideas. The main difference is
that instead of using the connection between the \AC\ polynomial and
$\gl$ he gives a combinatorial proof of Proposition~\ref{prop:gl}
similar to the proof of our lemmas in Section~4 of~\cite{FKV}.
I would like to thank him
for sending me a preliminary version of his preprint~\cite{Ch}.

I am thankful to D.~Bar-Natan, S.~Chmutov, 
J.~Figueroa-O'Farrill, T.~Kimura, and A.~Kricker
for useful comments.

\end{document}